\documentclass[english]{article}
\usepackage{amscd,amssymb,stmaryrd}
\usepackage[fleqn]{amsmath}
\usepackage{graphicx}
\usepackage{caption}

\usepackage[utf8]{inputenc}
\usepackage[export]{adjustbox}
\usepackage{wrapfig}
\usepackage{amstext}
\usepackage{pstricks,pst-node,pst-plot,pst-coil}
\usepackage{amsthm}
\usepackage{amsmath}
\usepackage{color}
\usepackage{mathtools}
\usepackage{hyperref}
\usepackage[english]{babel}
\usepackage{booktabs}
\usepackage{mathrsfs}
\usepackage{multirow}
\usepackage{subfiles}
\usepackage{subfigure}
\usepackage{lineno}
\usepackage{floatrow}
\usepackage{authblk}

\newfloatcommand{capbtabbox}{table}[][\FBwidth]

\providecommand{\keywords}[1]
{
  \small	
  \textbf{\textit{Keywords---}} #1
}
\title{Adaptive multiscale model reduction for nonlinear parabolic equations using GMsFEM}
 \author[1]{Yiran Wang}
 \author[1]{Eric Chung}
\author[2]{Shubin Fu}
\affil[1]{Department of Mathematics, The Chinese University of Hong Kong, Hong Kong SAR}
\affil[2]{Department of Mathematics, University of Wisconsin-Madison,WI, USA}
\begin{document}
\maketitle
%% -- ---------------------------------------------------------------------
%% -- ---------------------------------------------------------------------
\begin{abstract}
In this paper, we propose a coupled Discrete Empirical Interpolation Method (DEIM) and Generalized Multiscale Finite element method (GMsFEM) to solve nonlinear parabolic equations with application to the Allen-Cahn equation. The Allen-Cahn equation is a model for nonlinear reaction-diffusion process. It is often used to model interface motion in time, e.g. phase separation in alloys. The GMsFEM allows solving multiscale problems at a reduced computational cost by constructing a reduced-order representation of the solution on a coarse grid. In \cite{ref_url1}, it was shown that the GMsFEM provides a flexible tool to solve multiscale problems by constructing appropriate snapshot, offline and online spaces. In this paper, we solve a time dependent problem, where online enrichment is used. The main contribution is comparing different online enrichment methods. More specifically, we compare uniform online enrichment and adaptive methods. We also compare two kinds of adaptive methods. Furthermore, we use DEIM, a dimension reduction method to reduce the complexity when we evaluate the nonlinear terms. Our results show that DEIM can approximate the nonlinear term without significantly increasing the error. Finally, we apply our proposed method to the Allen Cahn equation.
%Last but not least, we give a brief review of Exponential Time Differencing (\textbf{ETD}) method, which can be applied to solve a particular parabolic equation: Allen-Cahn equation.

\end{abstract}
%% -- ---------------------------------------------------------------------
%% -- ---------------------------------------------------------------------
%\tableofcontents

\section{Introduction}

%\subsection{Objections}
In this paper, we consider the Generalized Multiscale Finite element method (GMsFEM) for solving nonlinear parabolic equations. The main objectives of the paper are the following: (1) to demonstrate the main concepts of GMsFEM and brief review of the techniques; (2) to compare various online enrichment techniques; (3) to discuss the use of the Discrete Empirical Interpolation Method (DEIM) and present its performance in reducing complexity.
%\subsection{The need for an adaptive model reduction approach}
GMsFEM is a flexible general framework that generalizes the Multiscale Finite Element Method (MsFEM) by systematically enriching the coarse spaces. The main idea of this enrichment is to add extra basis functions that are needed to reduce the error substantially. Once the offline space is derived, it stays fixed and unchanged in the online stage. In \cite{adaptive}, it is shown that a good approximation from the reduced model can be expected only if the offline information is a good representation of the problem. For time dependent problems, online enrichment is necessary. We compare two kinds of online enrichment methods: uniform and adaptive enrichment, where the latter focuses on where to add online basis. We will discuss it in numerical results with more details.
When a general nonlinearity is present, the cost to evaluate the projected nonlinear function still depends on the dimension of the original system,
resulting in simulation times that can hardly improve over the original system. One approach to reduce computational cost is the POD-Galerkin method \cite{ref8,ref9,ref10,ref11},
which is applied to many applications, for example, in \cite{ref12,ref13,ref14,ref15,ref16}.
DEIM focuses on approximating each nonlinear function so that
a certain coefficient matrix can be precomputed and, as a result, the complexity in
evaluating the nonlinear term becomes proportional to the small number of selected
spatial indices. In this paper, we will compare various approximations of the DEIM projection.
We will illustrate these concepts by applying our proposed method to the Allen Cahn equation.
The remainder of the paper is organized as follows. In section 2, we present the problem setting and main ingredients of GMsFEM. In section 3, we consider the methods to solve the Allen-Cahn equation.

\section{Multiscale model reduction using the GMsFEM}\label{sec:GMsFEM}
In this section, we will give the construction of our GMsFEM for nonlinear parabolic equations. First, we present some basic
notations and the coarse grid formulation in Section 2.1. Then, we present the construction of the multiscale snapshot functions and basis functions in Section 2.2. The online enrichment process is introduced in Section 2.3. %and an efficient way to handle the nonlinear term is illustrated in last subsection.\\

\subsection{Preliminaries}
Consider the following parabolic equation in the domain $\Omega \subset \mathbb{R}^d$
\begin{eqnarray}
\begin{aligned}
	\dfrac{\partial u}{\partial t} - \text{div}( \kappa \nabla u )  = & f& \quad \text{in }
	\Omega \times [0,T], \\
	u(x,0)  = &g(x)&  \quad \text{in } \Omega,\\
	u(x,t)  = &0& \quad \text{on } \partial \Omega \times [0,T].
\end{aligned}
\label{eqn:model}
\end{eqnarray}

Here, we denote the exact solution of (\ref{eqn:model}) by $u$, $\kappa(x)$ is a high-contrast and heterogeneous permeability field, $f = f(x,u)$ is the nonlinear source function depending on the $u$ variable, $g(x)$ is a given function and $T>0$ is the final time. We denote the solution and the source term at $t=t_n$ by $u(\cdot, t_n)$ and $f(u(\cdot,t_n))$ respectively.
 The variational formulation for the problem (\ref{eqn:model}) is: find $u(\cdot,t) \in H^1_0(\Omega)$ such that
 \begin{eqnarray}
\begin{aligned}
\left\langle \dfrac{\partial u}{\partial t}, v\right\rangle+\mathcal{A}(u,v)&=\left\langle f, v\right\rangle  \quad \text{in } \Omega\times [0,T], \quad\forall v \in H_0^1(\Omega),\\
		u(x,0)  &= g(x) \quad \text{in } \Omega,\\
	u(x,t)  &= 0 \quad \text{on } \partial \Omega \times [0,T].
\end{aligned}
\label{eqn:model_va}
\end{eqnarray}
where $\mathcal{A}(u,v)=\int_{\Omega} \kappa \nabla u\cdot \nabla v \;dx$.

  In order to discretize (\ref{eqn:model_va}) in time,
 we need to apply some time differencing methods. For simplicity, we first apply the implicit Euler scheme with time step $\Delta t>0$ and in section 3, we will consider a different differencing method: ETD. We obtain the following
 discretization for each time $t_n=n\Delta t,n=1,2,\cdots, N$ ($T=N\Delta t$),
 \begin{eqnarray*}
     \cfrac{u(\cdot,t_n)-u(\cdot,t_{n-1})}{\Delta t}=\text{div}(\kappa \nabla u(\cdot,t_n))+f(u(\cdot,t_n)).
 \end{eqnarray*}

\begin{figure}[ht]
    \centering
    \subfigure
       { \includegraphics[width=3in]{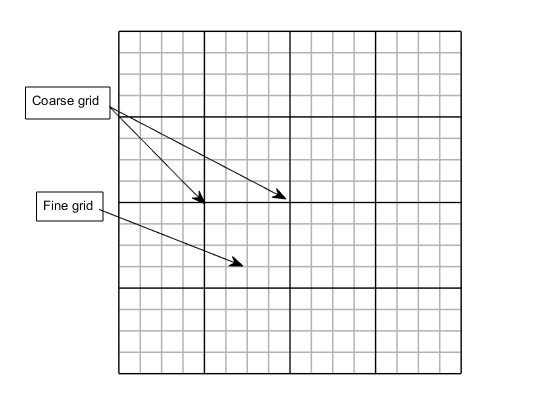}}
        \subfigure
       { \includegraphics[width=3in]{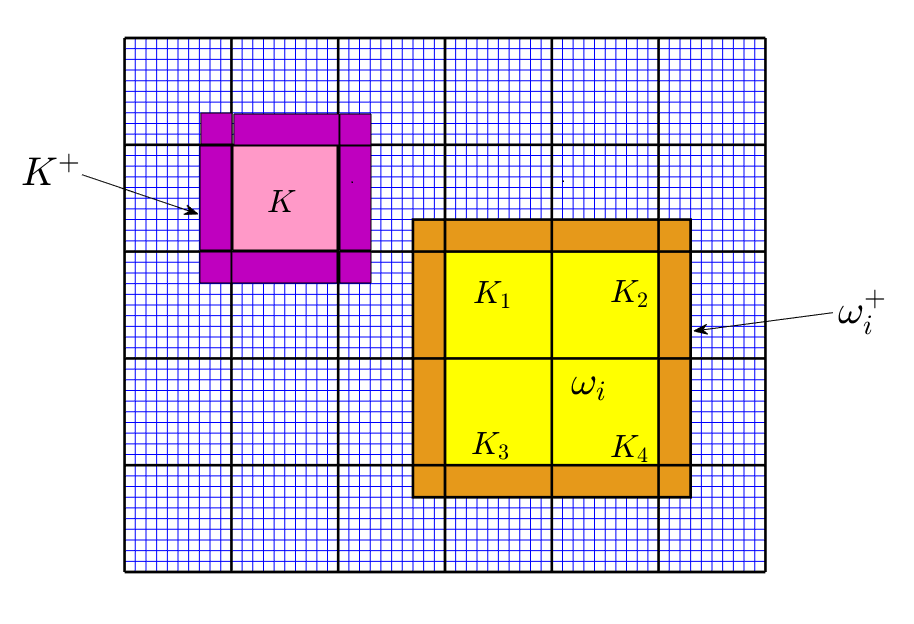}}

    \caption{Left: an illustration of fine and coarse grids. Right: an illustration of a coarse neighborhood, coarse element, and oversampled domain}\label{figure:coarse}
\end{figure}

Let $T^{h}$ be a partition of the domain $\Omega$ into fine finite elements. Here $h>0$ is the fine grid mesh size.
 The coarse partition, $T^{H}$ of the domain $\Omega$, is formed such that each element in $T^{H}$ is a connected union of fine-grid blocks. More precisely, $\forall K_{j} \in T^{H}$, $ K_{j}=\bigcup_{F\in I_{j} }F$ for some $I_{j}\subset T^{h}$. The quantity $H>0$ is the coarse mesh size. We will consider the rectangular coarse elements and the methodology can be used with general coarse elements. An illustration of the mesh notations is shown in the Figure \ref{figure:coarse}. We denote the interior nodes of $T^{H}$ by $x_i,i=1,\cdots,N_{\text{in}}$,
 where $N_\text{in}$ is the number of interior nodes. The coarse elements
 of $T^{H}$ are denoted by $K_j,j=1,2,\cdots,N_e$, where $N_e$ is the number of coarse elements. We define the coarse neighborhood of the nodes $x_i$ by $D_i:=\cup\{K_j\in T_{H}:x_i\in \overline{K_j}\}$.
\subsection{The GMsFEM and the multiscale basis functions}
 In this paper, we will apply the GMsFEM to solve nonlinear parabolic equations. The
 method is motivated by the finite element framework. First, a variational formulation is defined. Then we construct some multiscale basis functions.
  Once the fine grid is given, we can compute the fine-grid solution. Let $\gamma_1,\cdots,\gamma_n$ be the standard finite element basis, and define $V_f=\text{span}\{\gamma_1,\cdots,\gamma_n\}$ to be the
 fine space.
 We obtained the fine solution denoted by $u_f^n$ at $t=t_n$ by solving
  \begin{eqnarray}
  \begin{aligned}
      \frac{1}{\Delta t}\left\langle u_f^{n}, v\right\rangle+\mathcal{A}\left(u_f^{n}, v\right)&=\left\langle \frac{1}{\Delta t} u_f^{n-1}+f(u_f^{n}), v\right\rangle, \quad \forall v \in V_f,\\
      u_f^0&=g_h,
  \end{aligned}
  \label{fine sol}
  \end{eqnarray}
   where $g_h$ is the $V_f$ based approximation of $g$.
 The construction of multiscale basis functions follows two general steps. First, we construct snapshot basis functions in order to build a set of possible modes of the solutions. In the second step, we construct multiscale basis functions with a suitable spectral problem defined in the snapshot space. We take the first few dominated eigenfunctions as basis functions. Using the multiscale basis functions, we obtain a reduced model. %to solve the corresponding variations in space.

More specifically, once the coarse and fine grids are given, one may construct the multiscale basis functions to approximate the solution
of (\ref{eqn:model_va}). To obtain the multiscale basis functions, we first define the snapshot space. For each coarse neighborhood $ D_{i}$, define $J_h( D_{i})$ as the set of the fine nodes of $T^{h}$ lying on $\partial D_{i}$ and denote
the its cardinality by $L_i \in \mathbb{N}^{+}$. For each fine-grid node $x_j \in J_h( D_{i})$, we define a fine-grid function $\delta_{j}^{h}$ on $J_h( D_{i})$ as $\delta_{j}^{h}(x_k)=\delta_{j,k}$. Here
 $\delta_{j,k}=1$ if $j=k$ and $\delta_{j,k}=0$ if $j\neq k$. For each $j=1,\cdots, L_i$, we define the snapshot basis functions $\psi_{j}^{(i)}$ ($j=1,\cdots,L_i$) as the solution of the following system\\
\begin{eqnarray}
\begin{aligned}
   -\text{div}\left(\kappa \nabla \psi_{j}^{(i)}\right) &=0 \quad \text { in } D_{i} \\
   \psi_{j}^{(i)} &=\delta_{j}^{h} \quad \text { on } \partial D_{i}.\label{snap_basis}
\end{aligned}
\end{eqnarray}
 The local snapshot space $V_{\text { snap }}^{(i)}$ corresponding to the coarse neighborhood $ D_{i}$  is defined as follows
  $V_{snap}^{(i)}:=$ \text{span}$\{\psi_{j}^{(i)}:j=1,\cdots,L_{i}\}$ and the snapshot space reads $V_{\text {snap}} :=\bigoplus_{i=1}^{N_{\text {in}}} V_{\text {snap}}^{(i)}$.

  In the second step, a dimension reduction is performed on $V_{\text {snap}}$.
  For each $i=1,\cdots, N_{\text {in}}$, we solve the following spectral problem:
  \begin{eqnarray}
      \int_{D_{i}} \kappa \nabla \phi_{j}^{(i)} \cdot \nabla v=\lambda_{j}^{(i)} \int_{D_{i}} \hat{\kappa} \phi_{j}^{(i)} v \quad \forall v \in V_{\text {snap}}^{(i)}, \quad j=1, \ldots, L_{i}\label{eigen}
  \end{eqnarray}
  where $\hat{\kappa} :=\kappa \sum_{i=1}^{N_{i n}} H^{2}\left|\nabla \chi_{i}\right|^{2}$ and $\{\chi_{i}\}_{i=1}^{N_{i n}}$ is a set of partition of unity that solves the following system:
  \begin{eqnarray*}
  \begin{array}
  {rlrl}{-\nabla \cdot\left(\kappa \nabla \chi_{i}\right)} & {=0} & {} & {\text { in } K \subset D_{i}} \\ {\chi_{i}} & {=p_{i}} & {} & {\text { on each } \partial K \text { with } K \subset D_{i}} \\ {\chi_{i}} & {=0} & {} & {\text { on } \partial D_{i}}
  \end{array}
  \end{eqnarray*}
  where $p_i$ is some polynomial functions and we can choose linear functions for simplicity.
  Assume that the eigenvalues obtained from (\ref{eigen}) are arranged in ascending order and we may use the first $1<l_i \leq L_{i}$ (with $l_{i} \in
  \mathbb{N}^{+}$) eigenfunctions (related to the smallest $l_i$ eigenvalues) to
  form the local multiscale space $V_{\text{off}}^{(i)}:=$ snap$\{\chi_{i}\phi_{j}^{(i)}:j=1,\cdots,L_{i}\}$. The mulitiscale space $V_{\text{off}}^{(i)}$ is the direct sum of the local mulitiscale spaces,namely
  $V_{\text {off}} :=\bigoplus_{i=1}^{N_{\text {in}}} V_{\text {off}}^{(i)}$.
  Once the multiscale space $V_{\text {off}}$ is constructed, we can find the
  GMsFEM solution $u_{\text{off}}^n$ at $t=t_n$ by solving the following equation
  \begin{eqnarray}
  \begin{aligned}
        \frac{1}{\Delta t}\left\langle u_{\mathrm{off}}^{n}, v\right\rangle+ \mathcal{A}\left(u_{\mathrm{off}}^{n}, v\right)&=\left\langle \frac{1}{\Delta t} u_{\mathrm{off}}^{n-1}+f(u_{\mathrm{off}}^{n}), v\right\rangle, \quad \\
         \langle u_{\mathrm{off}}^{0},v\rangle&=\langle g,v\rangle, \quad\forall v \in V_{\mathrm{off}}.\label{ms sol}
  \end{aligned}
  \end{eqnarray}

\subsection{Online enrichment}
We will present the constructions of online basis functions \cite{ref1} in this section.
%The first way is adding same number of basis in online process, which is called uniform enrichment; The other one is adding basis only in the place where the corresponding residual exceeds some given tolerance, which is adaptive enrichment. We will present how to compute the residual as follows.
\subsubsection{Online adaptive algorithm}
In this subsection, we will introduce the method of online enrichment. After obtaining the multiscale space $V_{\text {off}}$, one may add some online basis functions based on local residuals.
Let $u_{\text {off}}^n \in V_{\text {off}}$ be the solution obtained in (\ref{ms sol}) at time $t=t_n$. Given a coarse neighborhood $D_i$, we define
$V_i:=H_0^1(D_i)\cap V_{\text {snap}}$ equipped with the norm
$\|v\|_{V_i}^{2}:=\int_{D_i}\kappa |\nabla {v}|^2$. We also define the local residual operator $R_i^n: V_i\rightarrow \mathbb{R}$ by
\begin{eqnarray}
    \mathcal{R}_{i}^{n}\left(v ; u_{\text{off}}^{n}\right) :=\int_{D_{i}}\left(\frac{1}{\Delta t} u_{\text{off}}^{n-1}+f(u_{\text{off}}^{n})\right) v-\int_{D_{i}}\left(\kappa \nabla u_{\mathrm{off}}^{n} \cdot \nabla v+\frac{1}{\Delta t} u_{\mathrm{off}}^{n} v\right), \quad \forall v \in V_{i}. \label{loc res}
\end{eqnarray}
 %where $u_{f}^{n-1}$ is the fine-scale solution at time
 %$t=t_{n-1}$.
 The operator norm $R_i^n$, denoted by $\|R_i^n\|_{V_{i}^{*}}$, gives a measure of the quantity of residual.  The online basis functions are computed during the time-marching process for a given fixed time $t=t_n$, contrary to
 the offline basis functions that are pre-computed.

  Suppose one needs to add one new online basis $\phi$ into the space $V_i$. The analysis in \cite{ref1} suggests that the required online basis $\phi\in V_i$ is the solution to the following equation
  \begin{eqnarray}
      \mathcal{A}(\phi, v)=\mathcal{R}_{i}^{n}\left(v ; u_{\text {off}}^{n, \tau}\right) \quad \forall v \in V_{i}.\label{online}
  \end{eqnarray}
  We refer to $\tau \in \mathbb{N}$  as the level of the enrichment
  and denote the solution of (\ref{ms sol}) by $u_{\text{off}}^{n,\tau}$.
  Remark that $V_{\text {off}}^{n,0}:=V_{\text{off}}$ for time level $n\in \mathbb{N}$. Let $\mathcal{I} \subset\left\{1,2, \ldots, N_{i n}\right\}$ be the index set over some non-lapping coarse neighborhoods. For each $i\in \mathcal{I}$, we obtain a online basis $\phi_i\in V_i$ by solving (\ref{online}) and define
$V_{\text {off}}^{n, \tau+1}=V_{\text {off}}^{n, \tau} \oplus \operatorname{span}\left\{\phi_{i} : i \in \mathcal{I}\right\}$.
After that, solve (\ref{ms sol}) in $V_{\text {off}}^{n, \tau+1}$.

%Consequently, following the arguments in (\cite{ref1}) ,we have at time $t=t_n$,
%\begin{eqnarray}
  %  \|u_f^n-u_{\text{off}}^{n,\tau+1}\|_V^2\leq(1-
   % \dfrac{\Lambda_{\text{min}}^{(\mathbb{I})}}{C_{\text{err}}}
    %\dfrac{\Sigma_{i\in\mathbb{I}}\|R_i^n\|_{V_i^{*}}(\lambda_{l_i+1}^{i})^{-1}}
   %{\Sigma_{i\in\mathbb{I}}\|R_i^n\|_{V_i^{*}}(\lambda_{l_i+1}^{i})^{-1}})
    %\|u_f^n-u_{\text{off}}^{n,\tau}\|_V^2 \label{on iter},
%\end{eqnarray}
%where $C_{\text{err}}$ is a uniform constant and $\Lambda_{\text{min}}^{(\mathbb{I})}=\text{min}_{i\in \mathbb{I}}\lambda_{l_i+1}^{i}$. Here, the norm is defined by
%$\|\cdot\|_{V}:=\sqrt{\mathcal{A}(\cdot,\cdot)}$. Inequality (\ref{on iter}) shows that we can obtain a better accuracy by adding more online basis functions at each time $t=t_n$ and the rate if convergence depends on the constant $C_{\text{err}}$ and
%$\Lambda_{\text{min}}^{(\mathbb{I})}$.

\subsubsection{Two online adaptive methods}
In this section, we compare two ways to obtain online basis functions which are denoted by online adaptive method 1 and online adaptive method 2 respectively.
%which will be illustrated in Table \ref{online algo1} and  Table \ref{online algo12}.
 Online adaptive method 1 is adding online basis using online adaptive method from offline space in each time step, which means basis functions obtained in last time step are not used in current time step. Online adaptive method 2 is keeping online basis functions in each time step. Using this accumulation strategy, we can skip online enrichment after a certain time period when the residual defined in (\ref{loc res}) is under given tolerance. We also presents the results of these two methods in Figure \ref{lalala} and Figure \ref{gegege} respectively.

\subsubsection{Numerical results}
\begin{figure}[ht]

\centering
\subfigure[$\kappa_1$]
{\includegraphics[width=3in]{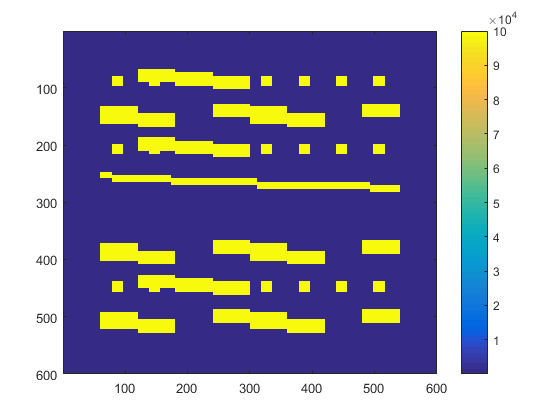}}
\subfigure[$\kappa_2$]
{\includegraphics[width=3in]{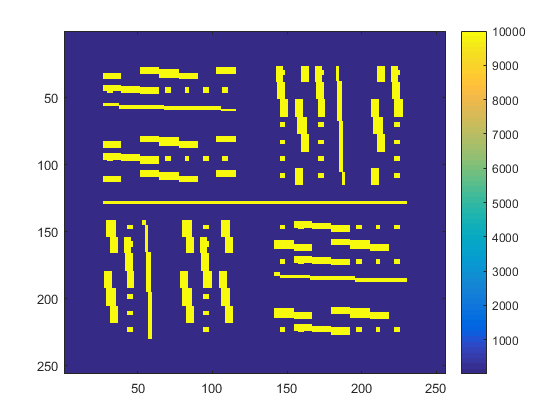}}
\caption{permeability field}\label{permeability}

\end{figure}

In this section, we present some numerical examples to demonstrate the efficiency of our proposed method. The computational
domain is $\Omega=(0,1)^2\subset \mathbb{R}^2$ and $T=1$. The medium $\kappa_1$ and $\kappa_2$ are shown in Figure \ref{permeability}, where the contrasts are $10^4$ and $10^5$ for $\kappa_1$ and $\kappa_2$ respectively. Without special descriptions, we use $\kappa_1$.

For each function to be approximated, we define the following quantities $e_a^n$ and $e_2^n$ at $t=t_n$ to measure energy error and $L^2$ error respectively.
 $$e_a^n=\frac{\left\|u_{\mathrm{f}}^{n}-u_{\mathrm{off}}^{n}\right\|_{V(\Omega)}}{\left\|u_{\mathrm{f}}^{\mathrm{n}}\right\|_{V(\Omega)}}\quad
  e_{2}^{n}=\frac{\left\|u_{\mathrm{f}}^{n}-u_{\mathrm{off}}^{n}\right\|_{L^{2}(\Omega)}}{\left\|u_{\mathrm{f}}^{\mathrm{n}}\right\|_{L^{2}(\mathcal{D})}}$$
  where $u_{\mathrm{f}}^{n}$ is the fine-grid solution (reference solution) and $u_{\mathrm{off}}^{n}$ is the approximation obtained by the GMsFEM method. We define the energy norm and $L^2$ norm of $u$ by
  $$\|u\|_{V(\Omega)}^2=\int_{\Omega}\|\nabla{u}\|^2
  \quad
  \|u\|_{L^2}^2=\int_{\Omega}\|u\|^2 .$$

\textbf{Example 2.1}.
 In this example, we compare the error using adaptive online method 1 and uniform enrichment under different numbers of initial basis functions.
 We set the mesh size to be $H=1/16$ and $h=1/256$. The time step is $\Delta t=10^{-3}$ and the final time is $T=1$. The initial condition is $u(x,y,t)|_{t=0}=4(0.5-x)(0.5-y)$. We set the permeability to be $\kappa_1$.
 We set the source term $f=\frac{1}{\epsilon^2}(u^3-u)$, where $\epsilon=0.01$. We present the numerical results for for the GMsFEM at time t=0.1 in Table \ref{tab_basis_1},\ref{tab_basis_2},and \ref{tab_basis_3}. For comparison, we present the results where online enrichment is not applied in Table \ref{tab_soourse}. We observe that the adaptive online enrichment converges faster.
Furthermore, as we compare Table \ref{tab_soourse} and Table \ref{tab_basis_1}, we note that the online enrichment does not improve the error if we only have one offline basis function per neighborhood. Because the first eigenvalue is small, the error decreases in the online iteration is small. In particular, for each iteration, the error decrease slightly. As we increase the number of initial offline basis, the convergence is very fast and one online iteration is sufficient to reduce the error significantly.

% Example 5.2  We keep H, h, time step, number of offline basis and the initial condition the same as in Example 5.1.  We set the source term to be $f=u^3+1$. In this time, we just \textbf{add one online basis each time step}, and result is shown below in Table \ref{onebasis}.\\

\begin{figure}[!htbp]
    \centering
    \subfigure[error with tolerance $10^{-4}$]
      { \includegraphics[width=3in]{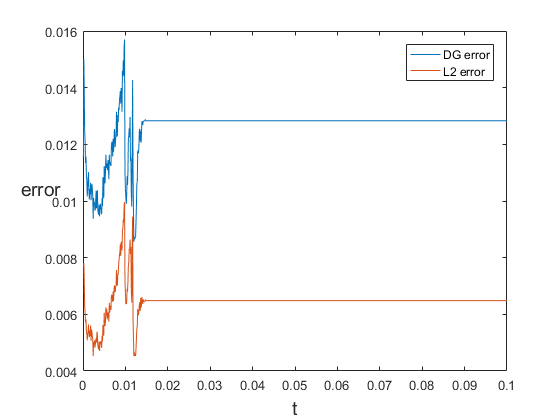}}
        \subfigure[DOF with tolerance $10^{-4}$]
      { \includegraphics[width=3in]{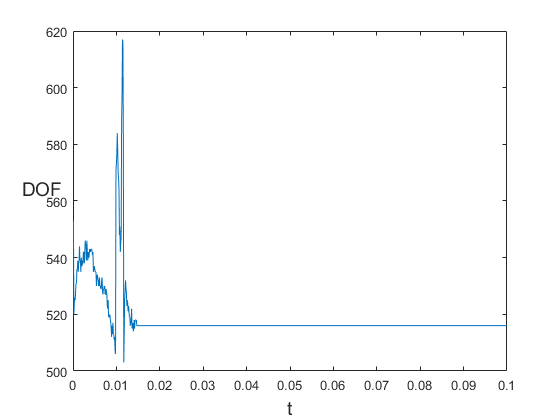}}
      	\quad
    \subfigure[error with tolerance $10^{-3}$]
	{\includegraphics[width=3in]{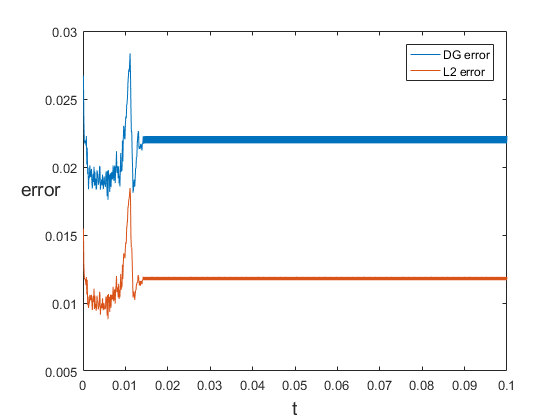}}
    \subfigure[DOF with tolerance $10^{-3}$]
	{\includegraphics[width=3in]{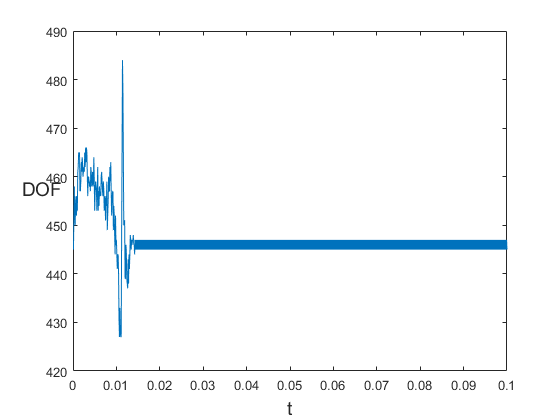}}
    \caption{error and DOF obtained by online method 1 in Example 2.2}\label{lalala}

\end{figure}

\begin{figure}[!htbp]
    \centering
    \subfigure[error with tolerance $10^{-4}$]
      { \includegraphics[width=3in]{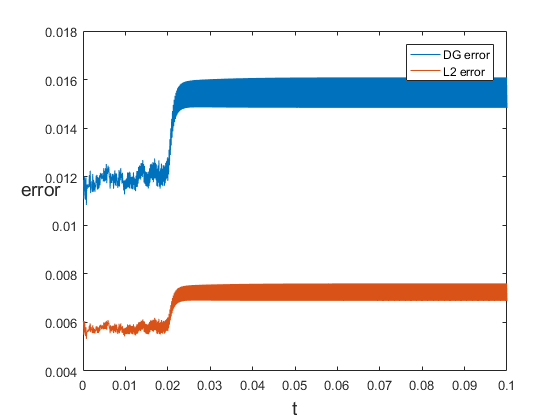}}
        \subfigure[DOF with tolerance $10^{-4}$]
      { \includegraphics[width=3in]{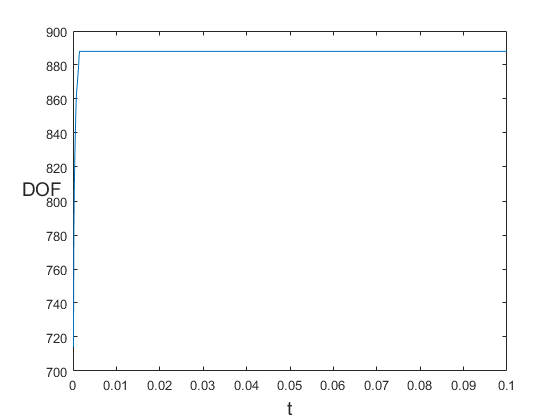}}

    \subfigure[error with tolerance $10^{-3}$]
	{\includegraphics[width=3in]{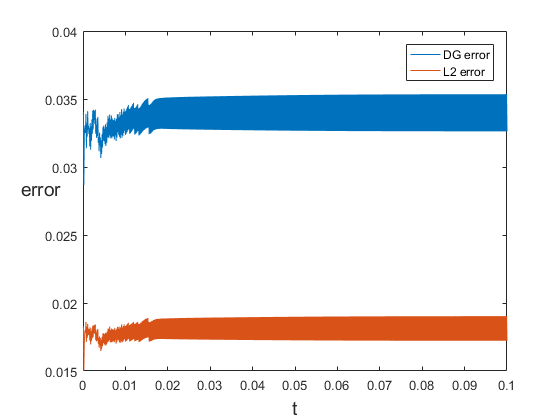}}
    \subfigure[DOF with tolerance $10^{-3}$]
	{\includegraphics[width=3in]{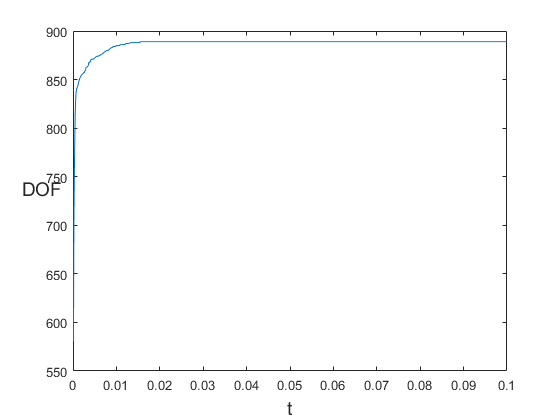}}
	\subfigure[error with tolerance $10^{-2}$]
	{\includegraphics[width=3in]{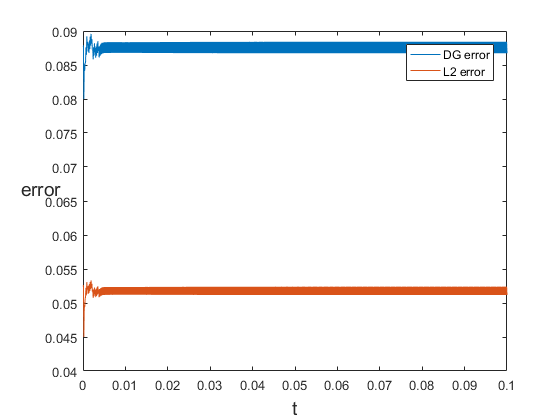}}
    \subfigure[DOF with tolerance $10^{-2}$]
	{\includegraphics[width=3in]{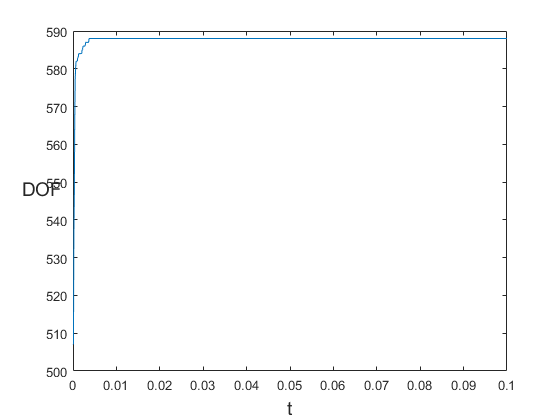}}
    \caption{error and DOF obtained by online method 2 in Example 2.2}\label{gegege}

\end{figure}

\begin{figure}[!htbp]
    \centering
    \subfigure[error with $\kappa_1$]
      { \includegraphics[width=3in]{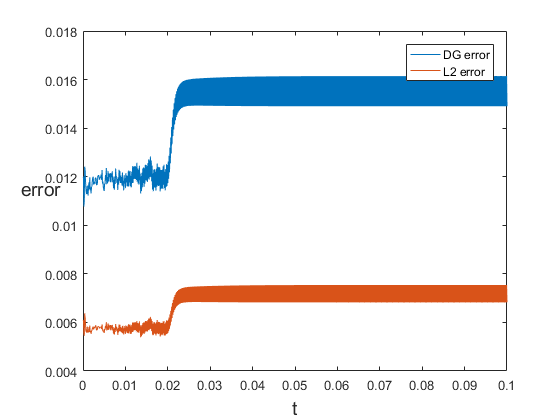}}
        \subfigure[DOF with  with $\kappa_1$]
      { \includegraphics[width=3in]{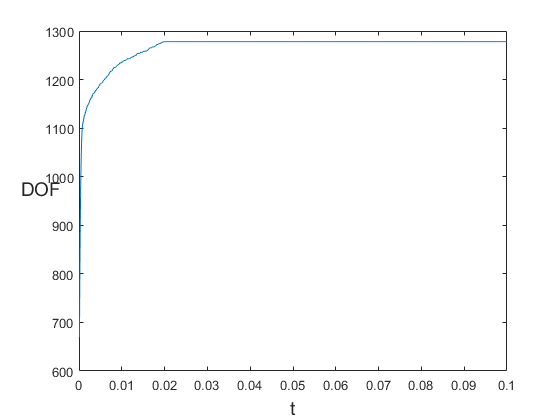}}
      \subfigure[error with $\kappa_2$]
      { \includegraphics[width=3in]{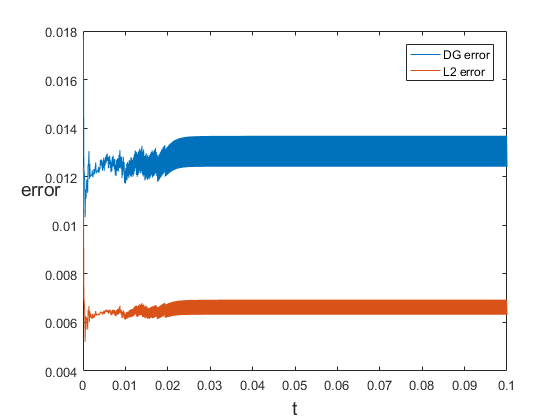}}
        \subfigure[DOF with  with $\kappa_2$]
      { \includegraphics[width=3in]{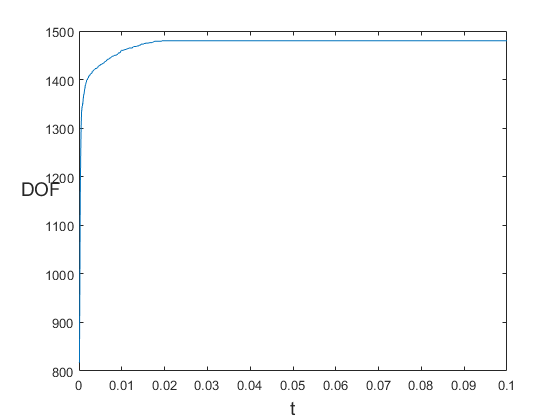}}
      \caption{error and DOF obtained by online method 2 in Example 2.2}\label{gee}
\end{figure}

\textbf{Example 2.2}.
We compare online Method 1 and 2 under different tolerance.
We keep H, h and the initial condition the same as in Example 2.1. We choose intial number of basis to be 450, which means we choose two initial basis per neighborhood.
We keep the source term as $f=\frac{1}{\epsilon^2}(u^3-u)$. When $\epsilon=0.01$, we choose the time step $\Delta t$ to be $10^{-4}$.
We plot the error and DOF from online Method 1 in Figure \ref{lalala} and compare with results from online Method 2 in Figure \ref{gegege}.
From Figure  \ref{lalala} and \ref{gegege}, we can see the error and DOF reached stability at $t=0.01$. In Figure \ref{gegege}, we can see the DOF keeps increasing before turning steady. The error remains at a relatively low level without adding online basis after some time. As a cost, online method 2 suffers bigger errors than method 1 with same tolerance.
We also apply our online adaptive method 2 under permeability $\kappa_2$ in Figure \ref{gee}. The errors are relatively low for two kinds of permeability.

\begin{table}[!htbp]
\caption{ The errors for online enrichment when number of initial basis=1. \\Left: Adaptive enrichment Right:Uniform enrichment.}
\begin{minipage}{0.48\linewidth}
\centering
	\begin{tabular}{c c c}
	\hline \hline
		DOF & $e_{a}$&$e_{2}$\\
	\hline
		225 & $14.47\%$&$19.55\%$ \\
		460 & $2.23\%$ & $1.14\%$\\
		550 & $1.20\%$ & $0.6\%$\\	
	
	\hline
	\end{tabular}
\end{minipage}
\begin{minipage}{0.48\linewidth}
	\begin{tabular}{c c c}
	\hline \hline
		DOF & $e_{a}$&$e_{2}$\\
	\hline
		225 & $14.48\%$&$19.56\%$ \\
		450 & $8.39\%$ & $6.54\%$\\
		675 & $2.45\%$ & $1.1\%$\\	
	
	\hline
	\end{tabular}
	\end{minipage}
    \label{tab_basis_1}
\end{table}

\begin{table}[!htbp]
\begin{minipage}{0.48\linewidth}
\centering
	\begin{tabular}{c c c}
	\hline \hline
		DOF & $e_{a}$&$e_{2}$\\
	\hline
		450 & $4.66\%$ & $2.64\%$\\
		681 & $1.65\%$ & $0.52\%$\\	
	\hline
	\end{tabular}
\end{minipage}
\begin{minipage}{0.48\linewidth}
	\begin{tabular}{c c c}
	\hline \hline
		DOF & $e_{a}$&$e_{2}$\\
	\hline
		450 & $4.65\%$ & $2.64\%$\\
		675 & $1.10\%$ & $0.669\%$\\	
	\hline
	\end{tabular}
	\end{minipage}
	\caption{The errors for online enrichment when number of initial basis=2. \\Left: Adaptive enrichment Right:Uniform enrichment}
	\label{tab_basis_2}
\end{table}

\begin{table}[!htbp]
\caption{The errors for online enrichment when  number of initial basis=3. \\Left: Adaptive enrichment Right: Uniform enrichment }
\begin{minipage}{0.48\linewidth}
\centering
	\begin{tabular}{c c c}
	\hline \hline
		DOF & $e_{a}$&$e_{2}$\\
	\hline
		675 & $2.89\%$ & $1.07\%$\\	
		903 & $0.944\%$ & $0.511\%$\\
	\hline
	\end{tabular}
\end{minipage}
\begin{minipage}{0.48\linewidth}
	\begin{tabular}{c c c}
	\hline \hline
		DOF & $e_{a}$&$e_{2}$\\
	\hline
		675 & $2.89\%$ & $1.07\%$\\	
		900& $1.13\%$ & $0.894\%$\\
	\hline
	\end{tabular}
	\end{minipage}
	\label{tab_basis_3}
\end{table}

\begin{table}[!htbp]
\caption{The errors for different $\epsilon$ in source term without online enrichment. \\ Up: energy error \quad Down: $L^2$ error }
\begin{minipage}{0.6\linewidth}
\centering
	\begin{tabular}{c|c|c}
	\hline \hline
		Source function&$ t=0.1$ & $t=0.2$ \\
	\hline
	 $\epsilon=0.1$&$5.97\%$&$5.94\%$\\
	 	$\epsilon=0.01$& $15.1\%$&$15.3\%$\\		
	
	\hline
	\end{tabular}
	\end{minipage}
	
\begin{minipage}{0.6\linewidth}
\centering
	\begin{tabular}{c|c|c}
	\hline \hline
		Source function&$ t=0.1$ & $t=0.2$\\
		 $\epsilon=0.1$& $4.57\%$&$4.57\%$\\	
	 $\epsilon=0.01$& $11.9\%$&$12.0\%$\\		
	\hline

	\hline
	\end{tabular}
		\end{minipage}
	\label{tab_soourse}
\end{table}

% \begin{table}[!htbp]

% \caption{The quantities adding one online basis function each time} \label{onebasis}
% 	\begin{tabular}{c|c|c|c}
% 	\hline \hline
% 		$ t=0.1$ & $t=0.2$ & $ t=0.4$& $t=0.8$\\
% 	\hline
% 	$1.642\%$&$1.649 \%$&$1.650\%$&$1.650\%$ \\

% 	\hline
% 	\end{tabular}

% \end{table}

\section{Application to the Allen-Cahn equation}
In this section, we apply our proposed method to the Allen-Cahn equation. We use the Exponential Time Differencing (ETD) for time dsicretization. To deal with the nonlinear term, DEIM is applied.  We will present the two methods in the following subsections.

\subsection{Derivation  of Exponential Time Differencing}

Let $\tau$ be the time step.
Using ETD, $u_{\text{off}}^{n}$ is the solution to (\ref{ms_etd})
     \begin{eqnarray}
\begin{aligned}
    \left\langle u_{\text{off}}^n,v\right\rangle+\tau\mathcal{A}(u_{\text{off}}^{n},v)&=\langle \text{exp}(-\dfrac{\tau}{\epsilon^2}\frac{f(u_{\mathrm{off}}^{n-1})}{u_{\mathrm{off}}^{n-1}})u_{\text{off}}^{n-1},v \rangle\\
  \langle u_{\mathrm{off}}^{0},v\rangle&=\langle g,v\rangle \quad\forall v \in V_{\mathrm{off}}\label{ms_etd}
\end{aligned}
\end{eqnarray}
Next, we will derive this equation.
%Since the term $\frac{1}{\epsilon^2}f(u)$ is stiff, if we solve it using explicit method, we need to take the time step $\tau$ to be $A\epsilon^2$ to make sure the stability of the scheme, where A is a constant. Although Backward Euler method is implicit, we can relax the strict restriction of time step $\tau$, but solving the nonlinear equation is difficult. To solve this problem, we consider exponential time differencing (ETD) method.\\
%The foundation of the method is actually integration factor method.
 %First we present the derivation of ETD method.
 %We choose a special source term $-\frac{1}{\epsilon^2}f(u(t_n,\cdot))$ in (\ref{eqn:model}), and we have
 We have
 \begin{eqnarray*}
u_t-\text{div}(\kappa \nabla u)+\frac{1}{\epsilon^2}f(u)=0
\end{eqnarray*}
  Multiplying the equation by integrating factor $e^{p(u)}$, we have
  \begin{eqnarray*}
e^{p(u)}u_t+e^{p(u)}\frac{1}{\epsilon^2}f(u)=e^{p(u)}\text{div}(\kappa \nabla u)
\end{eqnarray*}
 We require the above to become
  \begin{eqnarray}
\dfrac{d(e^{p(u)}u)}{dt}=e^{p(u)}\text{div}(\kappa \nabla u) \label{ode}
\end{eqnarray}
By solving
$$\dfrac{d(e^{p(u)}u)}{dt}=e^{p(u)}u_t+e^{p(u)}(\dfrac{d}{dt}p(u))u,$$
we have
$$p(u(t_n,\cdot))-u(0,\cdot))=\int_{0}^{t_n}\frac{1}{\epsilon^2}\frac{f(u)}{u}.$$
Using Backward Euler method in (\ref{ode}), we have
\begin{eqnarray}
        u_n-\tau \text{div}(\kappa \nabla u_n)=e^{-p(u)_{n}}u_{n-1}\label{half}
\end{eqnarray}
where $p(u)_n=p(u(t_n)-u(t_{n-1})).$
To solve (\ref{half}), we approximate (\ref{half}) as follows:
\begin{eqnarray}
    e^{-p(u)_{n}}u_{n-1}\approx e^{-\frac{\tau}{\epsilon^2}\frac{f(u(t_{n-1}))}{u(t_{n-1})}}u(t_{n-1}).
\end{eqnarray}
Using above approximation, we have
\begin{eqnarray}
         u_{\text{off}}^n-\tau \text{div}(\kappa \nabla  u_{\text{off}}^n)=\text{exp}(-\dfrac{\tau}{\epsilon^2}\frac{f(u_{\mathrm{off}}^{n-1})}{u_{\mathrm{off}}^{n-1}}) u_{\mathrm{off}}^{n-1}.\label{half_2}
\end{eqnarray}
%If we have $\frac{f(u)}{u}$ to be a constant term, then the approximation is exact. The advantage of ETD is that we can solve the linear term exactly.\\
%Similarly, we use GMsFEM to solve (\ref{half_2}) and we can obtain (\ref{ms_etd}).

  \subsection{DEIM method}
When we evaluate the nonlinear term, the complexity is $ O(\alpha(n)+c\cdot n)$, where $\alpha$ is some function and c is a constant. To reduce the complexity, we approximate local and global nonlinear functions with the Discrete Empirical Interpolation Method (DEIM)\cite{ref2}. DEIM is based on approximating a nonlinear function by means of an interpolatory projection of a few selected snapshots of the function. The idea is to represent a function over the domain while using empirical snapshots and information at some locations (or components).The key to complexity reduction is to replace the orthogonal projection of POD with the interpolation projection of DEIM in the same POD basis.

We briefly review the DEIM. Let $f(\tau)$ be the nonlinear function. We are desired to find an approximation of $f(\tau)$ at a reduced cost. To obtain a reduced order approximation of $f(\tau)$,
we first define a reduced dimentional space for it.
We would like to find m basis vectors (where m is much smaller than n), $\phi_1,\cdots,\phi_m$, such that we can write
$$f(\tau)=\Phi d(\tau),$$where $\Phi=(\phi_1,\cdots,\phi_m)$.
We employ POD to obtain $\Phi$ and use DEIM (refer Table \ref{DEIM algo}) to compute $d(\tau)$ as follows.
In particular, we solve $d(\tau)$ by using m rows of $\Phi$. This can be formalized using the matrix P
$$\mathrm{P}=\left[e_{\wp_{1}}, \ldots, e_{\wp_{m}}\right] \in \mathbb{R}^{n \times m},$$
where $e_{\wp_{i}}=[0,\cdots,1,0,\cdots,0]\in \mathbb{R}^{n}$ is the $\wp_i^{th}$ column of the identity matrix $I_n \in  \mathbb{R}^{n \times n}$ for $i=1,\cdots,m$.
Using $P^Tf(\tau)=P^T\Phi d(\tau)$, we can get the approximation for $f(\tau)$ as follows:
$$f(\tau) \approx \tilde{f}(\tau)=\Phi d(\tau)=\Phi\left(\mathrm{P}^{T} \Phi\right)^{-1} \mathrm{P}^{T} f(\tau)$$

\begin{table}[]

    \begin{tabular}{c l}
    \hline
       DEIM  & Algorithm  \\
       \hline
    \textbf{Input} & $\Phi=(\phi_1,\cdots,\phi_m)$ obtained by applying POD\\ &on a sequence of
                     $n_s$ functions evaluations\\
      \textbf{Output}&
      The interpolation indices $\overrightarrow{\lambda}=(\lambda_1,\cdots,\lambda_m)^T$\\
      &1. Set $[\rho,\lambda_1]=\max\{|\phi_1|\}$.\\
      &2. Set $\Phi=[\phi_1]$, $P=[e_{\lambda_1}]$, and  $\overrightarrow{\lambda}=(\lambda_1)$\\
      &3. for $i=2,\cdots,m$, do\\
      &Solve $(P^T\Phi)w=P^T\phi_i$ for some i.\\
      &Compute $r=\phi_i-\Phi w$\\
      &Compute $[\rho,\lambda_i]=\max\{|r|\}$\\
      &Set $\Phi=[\Phi,\phi_i]$, $P=[P,e_{\lambda_i}]$, and
$\overrightarrow{\lambda}=\left(\begin{array}{c}{\overrightarrow{\lambda}} \\ {\lambda_{i}}\end{array}\right)$\\
      & end for\\
      \hline
      \caption{DEIM algorithm}
      \label{DEIM algo}
      \end{tabular}
\end{table}
\subsection{Numerical results}
\begin{figure}[!htbp]
\centering
\subfigure[fine solution when $\epsilon=0.01$]
{\includegraphics[width=3in]{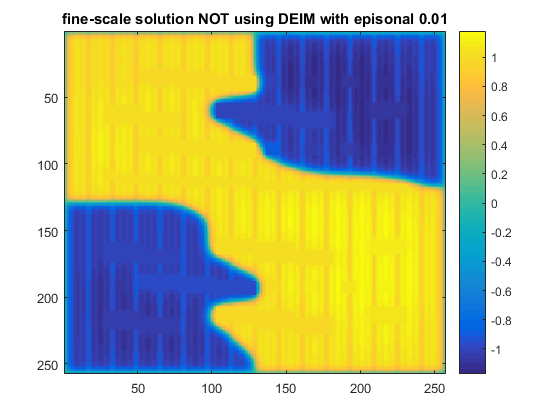}}
\subfigure[multiscale solution when $\epsilon=0.01$]
{\includegraphics[width=3in]{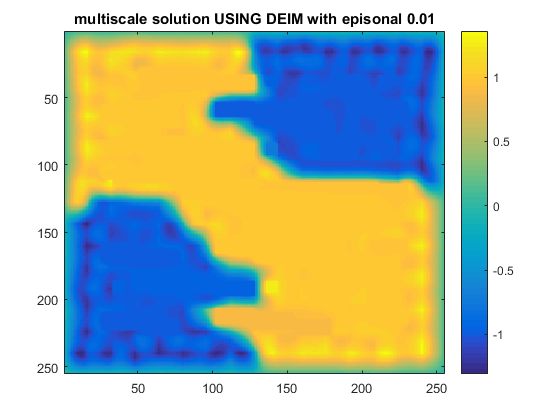}}
\caption{Comparing fine and multiscale solutions.}\label{haha}
\end{figure}

\textbf{Example 3.3}.
In this example, we apply the DEIM under the same setting as in Example 2.2 and we did not use the online enrichment procedure. We compare the results in Figure \ref{haha}.
%In terms of the construction of the snapshot $\Phi$, we use first solve a "similar" equation and combine the nonlinear functions in some time steps to get the snapshot. The "similar" equation differs from the original one in the following three cases.
To test the DEIM, we first consider the solution using DEIM where the snapshot are obtained  by the same equation.
First, we set $\epsilon=0.01$. We first solve the same equation and obtain the snapshot $\Phi$. Secondly, we use DEIM to solve the equation again. The two results are presents in Figure \ref{same}. The first picture are the errors we get when DEIM are not used while used in second one. The errors of these two cases differs a little since the snapshot obtained in the same equation.
Then we consider the cases where the snapshots are obtained:
\begin{enumerate}
    \item Different right hand side functions.
    \item Different initial conditions.
    \item Different permeability field.
    \item Different time steps.
\end{enumerate}
\subsubsection{Different right hand side}
\begin{figure}[!hp]
\centering
\subfigure[not using DEIM]
{\includegraphics[width=3in]{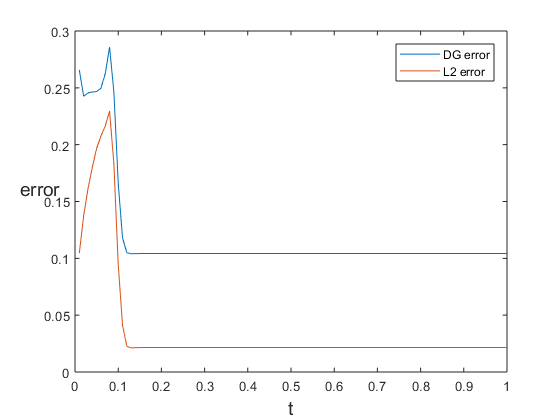}}
\subfigure[using DEIM]
{\includegraphics[width=3in]{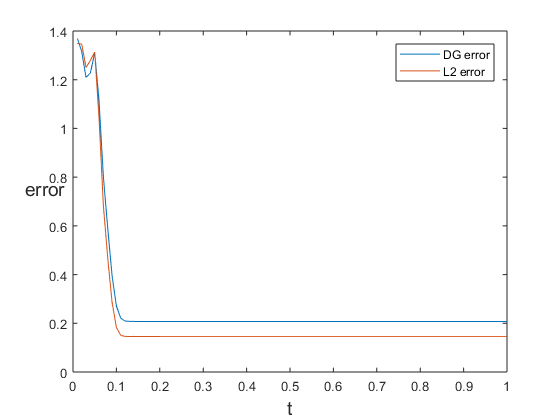}}
\caption{error for different $\epsilon$}\label{not same}
\end{figure}

Since the solution for different $\epsilon$ can have some similarities, we can use the solution from one to solve the other. In particular, since it will be more time-consuming to solve the case when $\epsilon$ is smaller. We can use the $f(u)$ for $\epsilon=0.09$ to compute the solution for $\epsilon=0.1$ since solutions for these two cases can only vary a little. I show the results in Figure \ref{not same}.

\subsubsection{Different initial conditions}
\begin{figure}[!htbp]

\centering
\subfigure[not using DEIM]
{\includegraphics[width=3in]{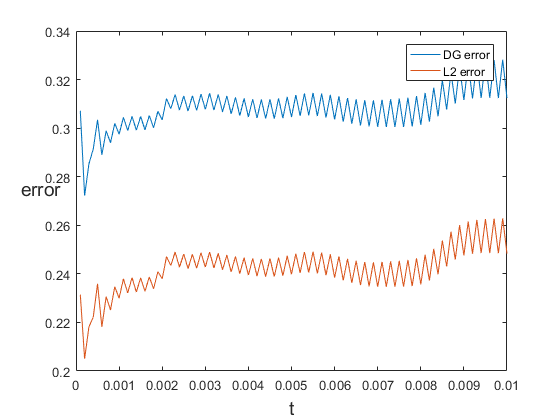}}
\subfigure[using DEIM]
{\includegraphics[width=3in]{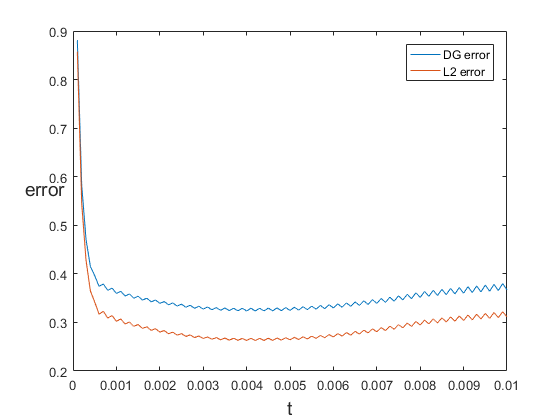}}
\caption{error for same $\epsilon$}\label{same}
\end{figure}

\begin{figure}[!htbp]
\subfigure[solution]
{\includegraphics[width=3in]{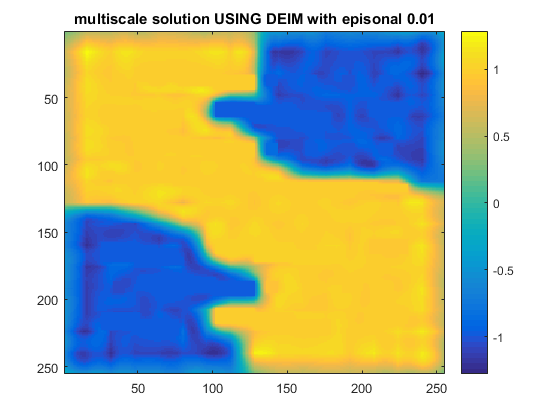}}
\subfigure[ error]
{\includegraphics[width=3in]{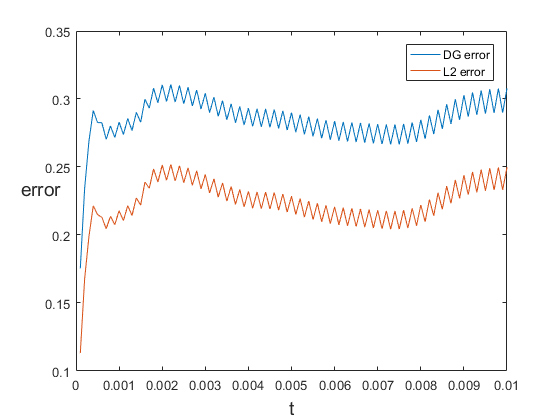}}
\caption{using DEIM for different initial conditions}\label{ini}
\end{figure}

In this section, we consider using the snapshot from different initial conditions, we record the results in Figure \ref{ini}. We first choose the initial condition to be compared Figure \ref{ini} and Figure \ref{same}, we can see that different initial conditions can have less impact on the final solution since the solution is close to the one where the snapshot is obtained in the same equation.
\subsubsection{Different permeability field }
\begin{figure}[!htbp]
\subfigure[fine solution under permeability $\kappa_2$]
{\includegraphics[width=3in]{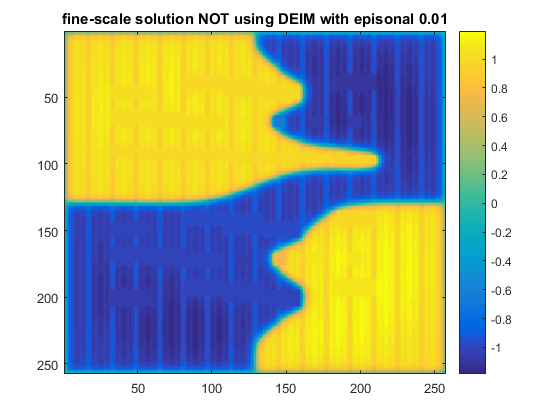}}
\subfigure[multiscale solution under permeability $a_2$]
{\includegraphics[width=3in]{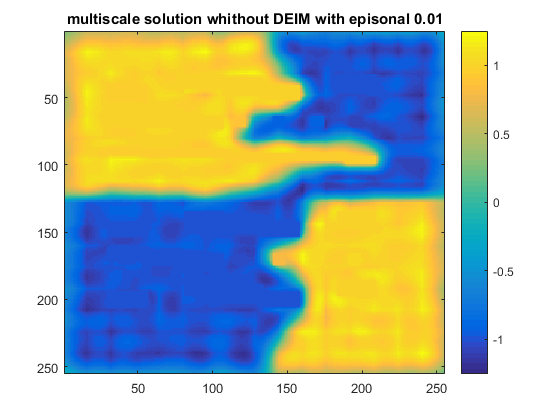}}
\subfigure[solution under permeability $\kappa_2$ using DEIM]
{\includegraphics[width=3in]{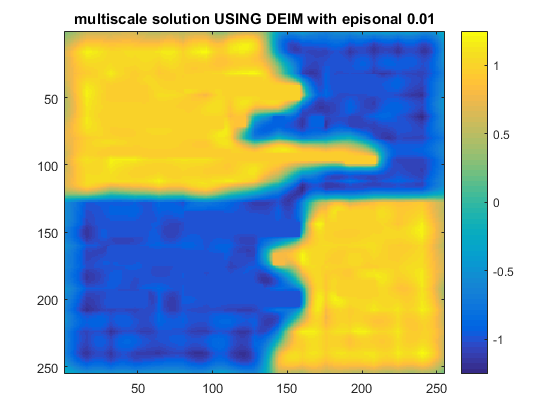}}
% \subfigure[error for multiscale solution]
% {\includegraphics[width=3in]{graphs/DEIM/error001_a2_ms.png}}
\subfigure[error for solution using DEIM]
{\includegraphics[width=3in]{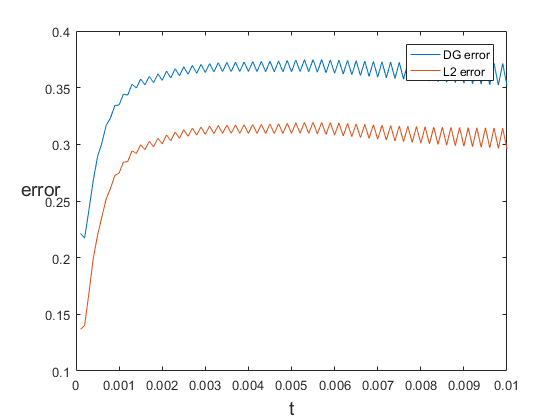}}
\caption{using DEIM for different permeability field}
\label{permeability_2}
\end{figure}

In this section, we consider using the snapshot from different permeability, we record the results in Figure \ref{permeability_2}. For reference, the first two figures plots the fine solution and multiscale solution without using DEIM. And we construct snapshot from another permeability $\kappa_1$ and we apply it to compute the solution in $\kappa_2$. The last figure shows the of using DEIM is relatively small.
\subsubsection{Different time steps}
\begin{figure}[!htbp]
\subfigure[solution]
{\includegraphics[width=3in]{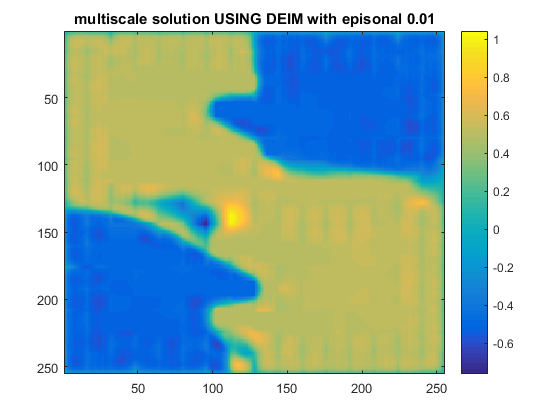}}
\subfigure[error]
{\includegraphics[width=3in]{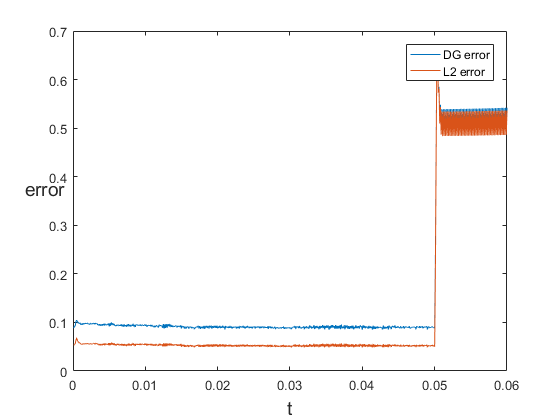}}
\caption{using DEIM for under different time step for $\kappa_1$}\label{step}
\end{figure}

\begin{figure}[!htbp]
\centering
\subfigure[solution]
{\includegraphics[width=3in]{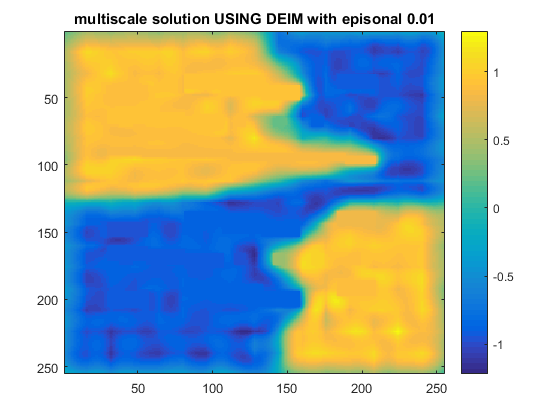}}
\subfigure[error]
{\includegraphics[width=3in]{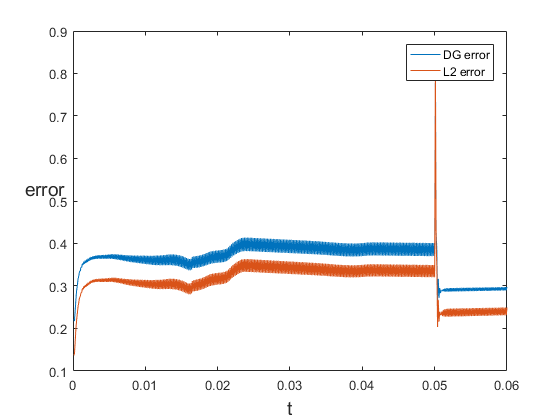}}
\caption{using DEIM for under different time step for $\kappa_2$\\}\label{step_1}
\end{figure}

In this section, we construct the snapshot by using nonlinear function obtained in previous time step for example when $t<0.05$. Then we apply it to DEIM to solve the equation in $0.05<t<0.06$. We use these way to solve the equation with permeability $\kappa_1$ and $\kappa_2$ respectively. We plot the results in Figure \ref{step} and \ref{step_1}. From these figures, we can see that DEIM have different effects applied to different permeability. With $\kappa_1$, the error increases significantly when DEIM are applied. But with  $\kappa_2$, the error decreased to a lower level when we use DEIM.


\section*{Acknowledgement}

The research of Eric Chung is partially supported by the Hong Kong RGC General Research Fund (Project numbers 14302018 and 14304719)
and CUHK Faculty of Science Direct Grant 2018-19.

% % \bibliographystyle{abbrv}
% % \bibliography{references}

\newpage
\begin{thebibliography}{8}
\bibitem{ref1}
Eric T.Chung, F., Yalchin Efendiev, S., Wing Tat Leung, T.: Residual-driven online generalized multiscale finite element methods.
 J. Comput. Phys \textbf{302}, 176–190(2015)

\bibitem{ref2}
Saifon Chaturantabut, F., Danny C. Sorensen, S.: Nonlinear model reduction via discrete empirical interpolation. SIAM J. Sci. Comput \textbf{32}(5),
2737–2764(2010)

\bibitem{adaptive}
Eric T. Chung, F.,Yalchin Efendiev, S., Thomas Y. Hou, T.: Adaptive multiscale model reduction with generalized multiscale finite element methods.
J. Comput. Phys \textbf{320}, 69-95(2016)
% \bibitem{ref3}
% T.Y. Hou and X.H.Wu, A multiscale finite element method for elliptic problems in composite
% materials and porous media, J. Comput. Phys., 134 (1997), pp. 169-189

\bibitem{ref4}
Yalchin Efendiev, F., Eduardo Gildin, S., Yanfang Yang, T.: Online Adaptive Local-Global Model Reduction for Flows in Heterogeneous Porous Media. Computation \textbf{4}, 22(2016)
% \bibitem{ref5}
% Thomas Y. Hou, Xiao-Hui Wu, and Zhiqiang Cai, Convergence of a multiscale
% finite element method for elliptic problems with rapidly oscillating
% coefficients, Mathematics of Computation, 68(227)(1999), pp.913-943
% \bibitem{ref6}
% Thomas Y. Hou, Xiao-Hui Wu, and Yu Zhang, Removing the cell resonance
% error in the multiscale finite element method via a Petrov-Galerkin
% formulation, Commun. Math. Sci., 2(2)(2004), pp.185-205
% \bibitem{ref7}
% Patrick Henning and Daniel Peterseim,Multiscale Model. Simul., 11(4), 1149–1175
\bibitem{ref8}
M. Hinze, F., S. Volkwein, S.:Proper orthogonal decomposition surrogate models for nonlinear dynamical systems: Error estimates and suboptimal control in Dimension Reduction of
Large-Scale Systems. Lect. Notes Comput. Sci. Eng. 45, pp. 261–306. Springer, Berlin(2005)
\bibitem{ref9}
K. Kunisch, F., S. Volkwein, S.: Control of the Burgers equation by a reduced-order approach using proper orthogonal decomposition. J. Optim. Theory Appl \textbf{102}, 345–371(1999)
\bibitem{ref10}
W. R. Graham, F., J. Peraire, S., K. Y. Tang, T.: Optimal control of vortex shedding using loworder models. Part I—Open-loop model development. Internat. J. Numer. Methods Engrg  \textbf{44}, 945–972(1999)
\bibitem{ref11}
S. S. Ravindran, F.: A reduced-order approach for optimal control of fluids using proper
orthogonal decomposition. Internat. J. Numer. Methods Fluids \textbf{34}, 425–448(2000)
\bibitem{ref12}
C. Prud’homme, F., D. V. Rovas, S., K. Veroy, T., L. Machiels, F., Y. Maday, F., A. T. Patera, S.,
G. Turinici, S.: Reliable real-time solution of parametrized partial differential equations:
Reduced-basis output bound methods. J. Fluids Eng  \textbf{124}, 70–80(2002)
\bibitem{ref13}
L. Machiels, F., Y. Maday, S., I. B. Oliveira, T., A. T. Patera, F., D. V. Rovas, F.: Output bounds
for reduced-basis approximations of symmetric positive definite eigenvalue problems. Comptes Rendus de l Académie des Sciences - Series I - Mathematics \textbf{331}, 153-158(2000)
\bibitem{ref14}
Y. Maday, F., A. T. Patera, S., G. Turinici, T.: A priori convergence theory for reduced-basis
approximations of single-parameter elliptic partial differential equations. J. Sci. Comput \textbf{17}, 437–446(2002)
\bibitem{ref15}
K. Veroy, F., D. V. Rovas, S., A. T. Patera, T.: A posterior error estimation for reduced-basis
approximation of parametrized elliptic coercive partial differential equations: “Convex
inverse” bound conditioners. ESAIM Control Optim. Calc. Var \textbf{8}, 1007–1028(2002)
\bibitem{ref16}
N. C. Nguyen, F., G. Rozza, S., A. T. Patera, T.: Reduced basis approximation and a posteriori
error estimation for the time-dependent viscous Burgers’ equation. Calcolo \textbf{46}, 157–185(2009)
\bibitem{ref17}
Gregory Beylkin, F., James M. Keiser, S., Lev Vozovoiy, T.: A New Class of Time Discretization Schemes for the Solution of Nonlinear PDEs. J. Comput. Phys \textbf{147}, 362–387(1998)
\bibitem{ref_url1}
Y. Efendiev, F., J. Galvis, S., T. Hou, T.: Generalized Multiscale Finite Element Methods. arXiv:1301.2866 [math.NA], \url{http://arxiv.org/submit/631572.}





\end{thebibliography}
% % \newpage
% % \section{Conclusion}

% % \begin{thebibliography}{}
% % \bibitem{ref1}
% % E. Chung, Y. Efendiev, W.T. Leung
% % \text{Residual-driven online generalized multiscale finite element methods}
% % J. Comput. Phys. 302 (2015) 176–190.\\

% % \bibitem{ref6}
% % \text{NONLINEAR MODEL REDUCTION VIA DISCRETE
% % EMPIRICAL INTERPOLATION}
% % SIAM J. SCI. COMPUT. Vol. 32, No. 5, pp. 2737–2764
% % \end{thebibliography}

%\section{Reference}

\end{document}